\documentclass[twocolumn,aps,secnumarabic,nofootinbib,floatfix]{revtex4}
\usepackage{times,graphicx}
\begin{document}
\title{Scholarly mathematical communication at a crossroads}
\thanks{Appeared as Nieuw Arch. Wisk. (5) \textbf{3} (2002), no. 3, 262--264}
\author{Greg Kuperberg}
\email[Email: ]{greg@math.ucdavis.edu}
\affiliation{Department of Mathematics,
    University of California, Davis, CA 95616}

\maketitle

Lately many librarians and some mathematicians have warned that academia faces
a serials crisis \cite{Kirby:journals}.  Ultimately I do not think that
scholarly mathematical communication is plagued by a crisis, but rather that it
is at a crossroads.  Computers and the Internet in general, and tools such as
TeX, MathSciNet, and the arXiv in particular, have enormously improved
mathematical communication. The question is whether, through leadership, we
will greatly extend these gains, or, through complacency, we will only see
marginal improvements. In my opinion, the math arXiv is a good foundation for
further progress.

Like Krzysztof Apt \cite{Apt:free}, I certainly think that the mathematical
research literature should be freely available.  But although subscription
prices are out of control, they are only one of several major shortcomings of
the traditional journal system.  The real impetus for change is that the old
system is slow and disorganized. All of the reform to date was first adopted by
non-paying researchers at privileged universities and institutes. So saving
money is ultimately only an important afterthought of the enterprise.

I especially agree with Apt that the key to progress is volunteer work and
public funds; they are what I mean by leadership.  But it is not so easy to use
these resources properly.  We won't succeed by simply moving the journal system
from the private sector to the public sector, and otherwise keeping it intact.
(So I partly agree with Michiel Kolman's skepticism on this point
\cite{Kolman:free}.) In this essay I will discuss three topics in the hope of
informing future decision-making:  The current state of the mathematical
literature, the math arXiv, and a proposal to reform peer review.

\section{The status quo}

Let us start with an idealization of the traditional system of scholarly
communication.  Mathematicians write uncertified, temporary preprints. They
submit them to journals, which edit them, bless them with peer review, and then
distribute them.  In fact this system was created more by circumstance than by
design, and it has always been supplemented by other venues of research
literature such as letters, books, and lecture notes.  Also, Zentralblatt and
Math Reviews have long organized an otherwise hopelessly disorganized array of
journals.

Electronic communication changed this picture in two important ways. First, it
made the informal preprint layer of the literature vastly more efficient. 
Second, and more importantly, it created enormous new pressure to further
organize the literature and accelerated a trend toward giantism.  As useful as
Math Reviews and Zentralblatt are on paper, they are far more useful on-line. 
The publishers are also agglomerating:  Elsevier (which now owns Academic
Press) and Springer both have large on-line libraries of mathematics papers
that other publishers may never be able to match.

The preprint trend and the giantism trend together also led to the arXiv, which
I will describe in more detail below.

First I would like to comment on the specter of giantism.   Naturally most
people don't like monopolies.  But working mathematicians want a universal
digital library of mathematics, and there is no clear way to sustain two of
them.  The same forces that propel Microsoft, America Online, and Amazon.com
are at play in this market.  So the question is not whether there should be a
centralized archive of the mathematical literature, but who should control it. 
I believe that a proprietary monopoly would be grossly against our interests,
even if it were owned by a non-profit professional society.  But with proper
oversight, a non-proprietary archive could remain useful and free.

One interesting precedent for such an archive is CTAN, the Comprehensive TeX
Archive Network \cite{ctan}.  All TeX users, including commercial publishers,
depend either directly or indirectly on CTAN.

\section{The arXiv}

The arXiv was created in 1991 by Paul Ginsparg as an e-mail-based service in
high-energy physics; it moved to the Web in 1993 \cite{arXiv}. (Ginsparg was
recently awarded a MacArthur Fellowship for the arXiv.) Now maintained at
Cornell, it was originally called hep-th, and it was for several years also
called the Los Alamos archives and the xxx e-print archive.  It does not claim
copyright to its articles, and all access is free.  It has some public
oversight in that it has advisory committees in each of its disciplines, and is
partly supported by the U.S. National Science Foundation. (I have heard
publishing figures criticize the public ``subsidy'' of the arXiv. For the
record, NSF award \#0132355 is only \$320,000 per year. This is a negligible
fraction of the tax dollars that are spent on math and physics journals.)

\begin{figure*}[ht]
\framebox{\scalebox{.75}{\includegraphics[clip=true,bb=24 288 588 768]
    {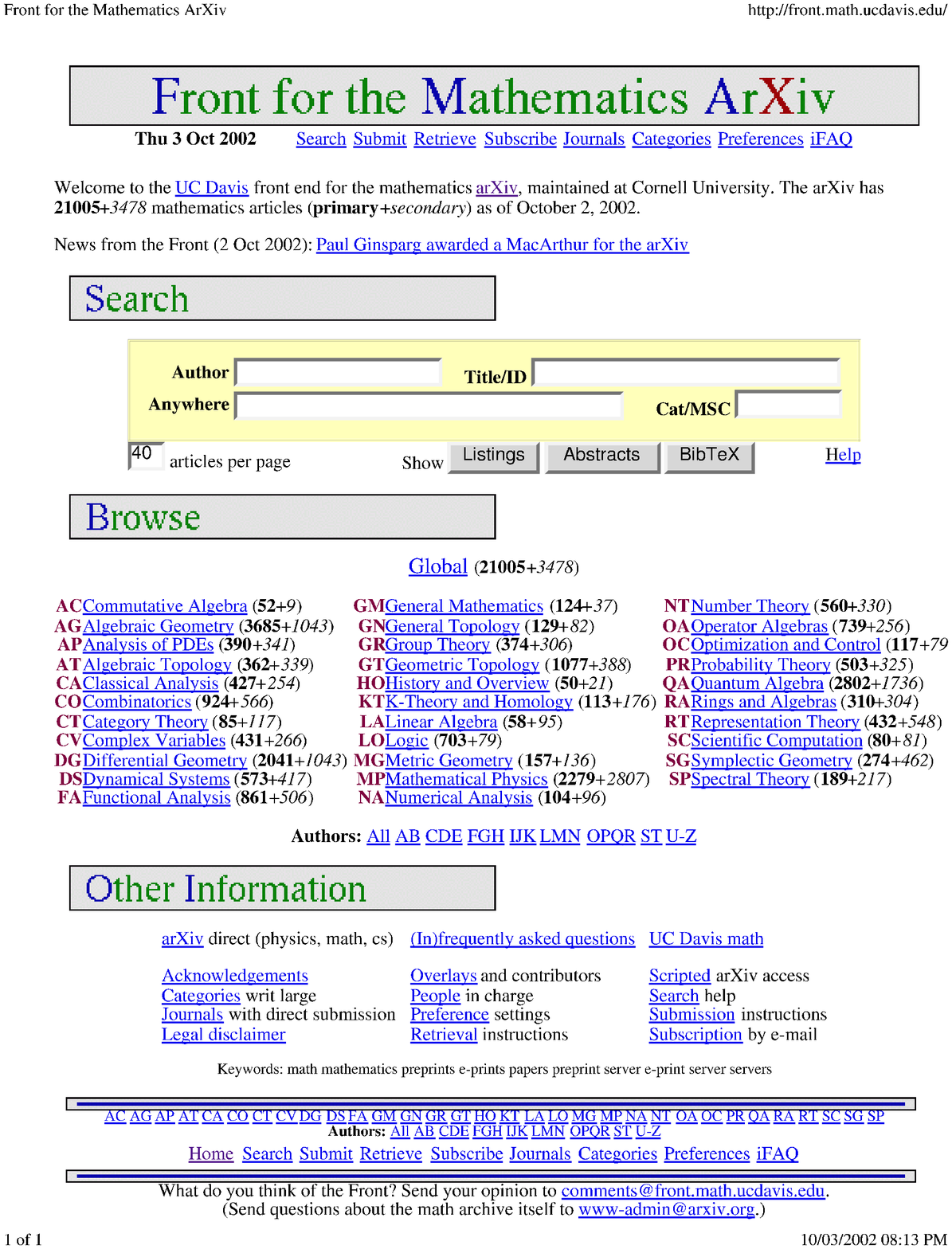}}}
\caption{The Front for the Mathematics ArXiv\cite{front}}
\label{f:front}
\end{figure*}

Authors usually submit articles to the arXiv directly, before submitting them
to journals. The arXiv distributes new submissions across the web and blesses
them as permanent at the end of each weekday.  Submissions may be revised, but
all versions of each article are numbered and old versions cannot be revoked.  
The arXiv is described in an article by Allyn Jackson \cite{Jackson:servers},
but readers who want to know more should also browse it on the web.  One
convenient means for doing so is the Front for the Mathematics ArXiv
(Figure~\ref{f:front}) \cite{front}.

The arXiv now has about 180,000 articles in physics, 20,000 in mathematics, and
nearly 3,000 in computer science. In many areas of physics, it is already a
universal library for the current literature.  Many people hope that it will
become just as important in mathematics.  The success of the math arXiv to date
is encouraging, but it will take many years for it to become universal.  As the
home page of the Front indicates, it is much more popular in some areas of
mathematics (algebraic geometry, quantum algebra, geometric topology) than in
others (applied mathematics, number theory). Finding ways to improve or further
popularize the arXiv will help us reach the ultimate goal sooner.

Significantly, the arXiv upsets the traditional dichotomy between unpublished
preprints and refereed publications.  Since arXiv articles are both permanent
and widely distributed, authors and readers take them just as seriously as
journal papers. (That is why arXiv documents are called articles, or e-prints,
and not preprints.)  For example, one of the most important papers in quantum
algebra is a still-unpublished arXiv article by Maxim Kontsevich
\cite{Kontsevich:poisson}.

A few journals are also arXiv overlays, which means that they submit their
published papers to the arXiv by proxy.  One of the most successful overlay
journals is \emph{Geometry and Topology} \cite{gandt}.  Here again the
distinction between preprints and publications blurs, since readers generally
do not keep track of authors' versus journals' versions of papers.

Nonetheless the arXiv does not completely displace journals. They are no longer
needed to typeset, distribute, or archive the mathematical literature, but they
still perform the crucial function of peer review.  This is why almost all
arXiv articles are still submitted to journals or conference proceedings.
Overlay journals also illustrate this point: While journal status does not
distinctly affect the text of an arXiv article, the fact that an article has
been accepted by a journal does still matter to readers, and even more to
authors.

Serious peer review probably cannot fit directly into the arXiv. The arXiv is
lightly moderated to screen out miscued and misclassified submissions (which is
why submissions are not announced immediately), but this is a superficial
effort compared to refereeing papers, or to reviewing them for Math Reviews or
Zentralblatt.  Rather, in heavily arXived areas of research, journals are a
\emph{de facto} second layer of the permanent literature, with the arXiv as the
first layer.  (Hence the term ``overlay''.)

\section{Peer review}

In the idealized journal system, the diligent referee first checks the main
results of a submitted paper.  If they are correct, the referee and the editor
then consider whether the results meet the journal's standard. In practice the
system is far from ideal.  Referees have very little accountability (although
some do an admirable job anyway).  Many papers are accepted on the basis of
name recognition or out of guilt -- just because the referees sat on them for
too long.  Authors need not take no for an answer, because they can scout for
journals that will publish them.  Journal editors should serve as a second line
of defense, but in practice they are only slightly more accountable than
referees.  (Again, some do an admirable job anyway.)  In the end, readers do
not know who refereed any given paper or why it was accepted.  Because papers
can only be published once, the system reduces peer review to simple binary
approval.

In my opinion, Math Reviews and Zentralblatt are inherently more useful forms
of peer review, because the reviews are not anonymous.  Their publishers do not
agree with me; by their own rules, reviewers are not required to fully referee
papers.  But some do so anyway.  Many mathematicians know of notoriously
mistaken papers that were inexplicably published.  The community is often angry
with the referees of such papers, but anonymity protects them from the readers
rather than the authors. Typically the Math Review sets the record straight. 
My favorite such case is the celebrated review by G\'abor Fejes T\'oth of Wu-Yi
Hsiang's inadequate proof of the Kepler sphere-packing conjecture
\cite{MR95g:52032}.  Even in the usual case when both authors and referees do a
good job, Math Reviews obviously inform readers more than any referee reports
could, since the latter are confidential.

In the presence of the arXiv, it is relatively easy to reform journals so that
they function much more like Math Reviews.  If a journal is purely an arXiv
overlay, then it need not take possession of its papers. So why should it wait
for authors to submit to it?  It could instead allow anyone to nominate
(``submit") any arXiv article for review, whether or not it has been published
elsewhere.  Let us call such a review service an ``open journal''.  If an open
journal reduces selection to its utilitarian minimum, it should add some other
value for readers to take it seriously.  It is therefore natural for open
journal referees to write public, non-anonymous reports like those in Math
Reviews.  (But reviewers who are not interested in a submission can reject it
privately and anonymously.)

Open journals have been tried before, both in connection with the arXiv and
elsewhere.  For example Quick Reviews \cite{quickrev} is a review service in
quantum computation maintained by Daniel Gottesman; it is mostly but not
exclusively arXiv-based.  But existing experiments lack a crucial feature: 
They are not designed to substitute for journal names in the author's list of
publications. For this purpose an open journal should do three things. It
should retain the trappings of a traditional journal, such as an editorial
board, an on-line masthead, and a bibliographic citation style.  It should also
keep the author informed of the status of papers under review.  And it should
prod editors and reviewers to attend to submissions, as traditional journals
do.

While many readers presume that referees check the results of papers, in
practice editors would scare away their referees if they actually demanded
this.  This inconsistency is only tenable because refereeing is anonymous. My
best idea to address the problem is to have the reviewers check one of three
options:
\begin{enumerate}\setlength{\itemsep}{0pt}
\item[\bf 1.] I have checked the main results.
\item[\bf 2.] I do not doubt the main results.
\item[\bf 3.] I doubt the main results.
\end{enumerate}
Some reviews might need to be co-signed to support option 1. Presumably option
3 would be rarely used.

It remains to be seen whether the open journal model will succeed.  Or one
might ask, who will pay people to do all of the work?  I see two reasons to be
cautiously optimistic.  First, as we already know, publishers pay the
mathematicians who manage, edit, and referee traditional journals very little
or not at all.  Most of the paid work for a traditional journal is for
typesetting and distribution. Most of the incentive for the mathematicians
involved is in credit for professional service.  In principle this incentive
could also sustain open journals. Second, I think that once the arXiv is
sufficiently established in mathematics, journals will naturally evolve toward
the open model.  But I hope that, through leadership, we can reach the future
sooner.

% \bibliography{soft,qa,mg}

\providecommand{\bysame}{\leavevmode\hbox to3em{\hrulefill}\thinspace}

\end{document}